\documentclass[12pt]{article}%
\usepackage[english]{babel}
\usepackage[utf8]{inputenc}
\usepackage[a4paper, top=2.5cm, bottom=2.5cm, left=2.2cm, right=2.2cm]%
{geometry}

\usepackage{times}
\usepackage{amsmath}
\usepackage{amsthm}
\usepackage{bbm}
\usepackage{changepage}
\usepackage{amssymb}
\usepackage[english]{babel}
\usepackage[utf8]{inputenc}
\usepackage{graphicx}

\newtheorem{theorem}{Theorem}

\newtheorem{lemma}[theorem]{Lemma}

\def\R{\mathbb{R}}
\def\C{\mathbb{C}}

\def\1{\mathbbm{1}}

\usepackage[
backend=biber,
style=alphabetic,
sorting=ynt
]{biblatex}
\begin{document}
\title{Cutoff Phenomenon and Limiting Profile of a Random Walk on the Symmetric Group}
\author{Ahmed Farah}
\date{May 12th 2021}
\maketitle
\section{Introduction}
When shuffling a deck of cards, we're often interested in in knowing the number of shuffles after which the deck is ``sufficiently shuffled," i.e. given a known starting configuration of the card deck and some probabilistic shuffling algorithm, how many times should we shuffle the deck so that the probability of any card being in any position is (almost) exactly the same? \\
We will consider this basic shuffling algorithm: Using a deck of $n$ cards, laid out and facing down, pick a card A uniformly at random. Pick another one, B, uniformly at random and swap their position. In probability parlance, this defines a probability measure on the set of possible permutations of the deck, which is the symmetric group $S_n$. It is defined as:
\begin{equation} \label{def:prob}
P(g)=\begin{cases}
\frac{1}{n}& \text{if $g$ is the identity}\\
\frac{2}{n^2}& \text{if $g$ is a transposition}\\
0 & \text{otherwise}
\end{cases}
\end{equation}
We can think of this process as a Markov chain, or a random walk on the Cayley graph of $S_n$, generated by the set of transposition. Hence, taking time steps in the random walk, i.e. repeating the shuffle, corresponds to the convolution of $P$ with itself, which we denote as:
\begin{align} \label{eq:1}
    P^{*2}(s) = (P * P)(s) = \sum_{t \in S_n}P(t^{-1})P(ts)
\end{align}
Let $P^{*k}$ denote the convolution of $P$ with itself $k$ times. The question of how far we are from a shuffled/random deck then turns into studying the behavior of the total variation norm as $k$ grows larger:
\begin{align} \label{eq:2}
||P^{*k} - U||_{TV} =: \underset{A \subseteq S_n}{\max} |P^{*k}(S) - U(S)| = \frac{1}{2}\sum_{s \in S_n} |P^{*k}(s) - \frac{1}{n!}|
\end{align}
Here, $U$ denotes the uniform distribution on $S_n$. \\
The main result that we will be discussing is an upper bound for this norm developed by Diaconis and Shahshahani \cite{DS}:

\begin{theorem}(Diaconis and Shahshahani, 1980): \\
Let $c$ be some real constant and $k = \lfloor \frac{1}{2} n \log n + cn \rfloor$. Assuming that $n \geq 10$, there exists a constant $a$ s.t.
\begin{align} \label{thm:DS}
||P^{*k} - U||_{TV} \leq ae^{-2c}
\end{align}
\end{theorem}
\noindent In other words, after $\frac{1}{2} n \log n$ convolutions, the distance decays exponentially. Here, as well as in the rest of the paper, we make the convention of using $\log$ to refer exclusively to the natural logarithm. Going forward, we will also omit floors when it is clear that the argument ought to be an integer, and use $|| \cdot ||$ to refer to the Total Variation distance, unless otherwise noted.\\
What's remarkable about Diaconis and Shahshahani's method is that it capitalizes on the rich representation theory of the symmetric group, turning the problem of asymptotically bounding the distance between two probability measures into one of bounding the character ratios of the representations of the symmetric group. This will be further explored in section 2.\\
\\
\noindent This theorem motivates another interesting question. Letting $f_n(c) = ||P^{*k} - U||$, we can define $f(c) = \underset{{n \to \infty}}{\lim} f_n(c)$ and question the behavior of $f$. From Theorem 1, we know that $f$ must decay exponentially in $k$. In 2020 Teyssier \cite{Teyssier} was able to determine $f$ exactly:
\begin{theorem}(Teyssier, 2020) \\
Let $c$, $k$, $P$ be defined as in Theorem 1. Let Poiss$(\lambda)$ denote a Poisson random variable with parameter $\lambda$, then:
\begin{align} \label{thm:T}
||P^{*k} - U||_{TV} \underset{n \to \infty}{\longrightarrow} ||\text{Poiss}(1+e^{-2c}) - \text{Poiss}(1)||_{TV}
\end{align}
\end{theorem}
\noindent In this paper, our goal is to present a high level overview of the proof techniques of both of these results, highlighting the significance of some intermediary results as we go along. Section 2 will center around proving the Upper Bound Lemma, a standard result which allows us to translate the problem of bounding the norm in (\ref{eq:2}) into the language of representation theory. Section 3 gives an overview of Diaconis and Shahshahani's upper bound on the variation distance in (\ref{thm:DS}). Section 4 gives an overview of some of the techniques used in Teyssier's proof in (\ref{thm:T}).

\section{The Upper Bound Lemma}
The goal of this section is to present a fairly self-contained proof of the upper bound lemma, which we state below.
\begin{theorem} \label{UBL}
Let $G$ be a finite group, $P$ be some probability distribution on $G$, and $U$ be the uniform distribution, then:
\begin{equation}
    ||P-U||^2 \leq \frac{1}{4} \sum_{\rho \in \widehat{G}^*} d_\rho \: \text{Tr}(\widehat{P}(\rho)\widehat{P}(\rho)^*)
\end{equation}
\end{theorem}
\noindent The notation and terminology used is explained in the following subsection. Serre Chapters 1 and 2 \cite{Serre} contain a more detailed overview.
\subsection{Representation theory preliminaries}
Let $G$ be a finite group. A (linear) \textbf{representation} $\rho$ is a homomorphism which maps $G$ to $GL_n(\C)$. Here, $n$ is referred to as the \textbf{degree} of $\rho$ and denoted $d_\rho$. The \textbf{character} of $\rho$, denoted as $\chi_\rho: G \to \C$ is a mapping defined by $\chi_\rho(g) =: \text{Trace}(\rho(g))$. More generally, let $\C[G]$ be the set of all functions $f: G \to \C$. Let $\widehat{G}$ denote the set of \textbf{irreducible representations} of $G$, and $\widehat{G}^*:= \widehat{G}/\{\text{triv}\}$, where \emph{triv} is the trivial representation \\
By the rotation invariance of the trace, we can see that the character satisfies a notable property: it is constant on the conjugacy classes of $G$. We call any such function a \textbf{class function}. Another notable example of a class function is the probability measure defined in (\ref{def:prob}): Remember that the conjugacy classes of $S_n$ are determined by cycle shapes, so the identity element is its own singleton conjugacy class, and the transposition together form a conjugacy class.

\subsection{The Fourier Transform}
Let $f$ be a function on $G$. We define the Fourier transform of $P$ at the representation $\rho$ as:
\begin{equation} \label{def:ft}
    \hat{f}(\rho) = \sum_{s \in G} f(s)\rho(s)
\end{equation}
What is notable is that when $f$ is a class function $G \to \C$, it commutes with $\rho$. To see this, fix $g \in G$ and note that:
\begin{align*}
    \rho(g) \hat{f}(\rho)
    & =  \rho(g) \sum_{s \in G} f(s)\rho(s) \\
    & =  \sum_{s \in G} f(s)\rho(gs) \\
    & =  \sum_{s \in G} f(s)\rho(gsg^{-1}) \rho(g) \\
    & =  \sum_{t \in G} f(t)\rho(t) \rho(g) \hspace{12pt}(\text{by setting $t = gsg^{-1}$ )} \\
    & =  \hat{f}(\rho) \rho(g) \\
\end{align*}
This means that by Schur's lemma that it's a homothety, i.e. $\hat{f}(\rho) = \lambda I$ for some $\lambda \in \C$ \\
The reason we want to study Fourier transforms of probability measures is because they turn convolution, a rather complicated operation, into multiplication, a simpler one. We mean this in the following sense:
\begin{lemma}
Let $f,h \in \C[G]$, then $\forall \rho$: $\widehat{f*h}(\rho) = \widehat{f}(\rho)\widehat{h}(\rho)$
\end{lemma}
\noindent \textit{Proof:}
\begin{align*}
    \widehat{f*h}(\rho)
    & = \sum_{s \in G} (f*h)(s) \rho (s)
    = \sum_{s \in G} \sum_{t \in G} f(t) h(t^{-1}s) \rho (s)
    = \sum_{s \in G} \sum_{t \in G} f(t) \rho (t) h(t^{-1}s) \rho (t^{-1}s) \\ 
    & = \sum_{t \in G} f(t) \rho (t) \sum_{s \in G} h(t^{-1}s) \rho (t^{-1}s)
    = \widehat{f}(\rho)\widehat{h}(\rho)
\end{align*}
Here, $f,h$ need not be class-invariant probability distributions, but any functions in $\C[G]$. When $P$ is a class function, however, then this statement combined with Schur's lemma implies that $\widehat{P^{*k}} = (\widehat{P})^k = \lambda^k I$. \\
Notice that $\forall f \in \C[G]$, $\hat{f}(\text{triv}) = \sum_{s\in G} f(s)$. In particular, when $f$ is a probability distribution, $\hat{f}(\text{triv}) = 1$, so of course, $\text{Tr}(\hat{f}(\text{triv})) = 1$,. \\
On the other hand, pick $\rho \in \widehat{G}$ and note that:
\begin{equation} \label{value:U}
\hat{U}(\rho) = 
\begin{cases}
1 & \text{when $\rho$ = triv} \\
0 & \text{otherwise}
\end{cases}
\end{equation}
Where $U$ is the uniform distribution. There is a slight abuse of notation in the statement, which is that the $0$ above refers to the $0$-operator or matrix on a $d_\rho$ dimensional space. The first case follows from what we've just shown. To see why the second holds, note that $\hat{U}(\rho) = \lambda I$ by Schur's lemma, and that by the orthogonality of the characters of irreducible representations, we get:
\[
\text{Tr}(\hat{U}(\rho)) = \lambda d_\rho = \frac{1}{|G|} \sum_{s \in G} \chi_\rho(s) = (\chi_\rho|\chi_\text{triv}) = 0
\]
\subsection{Plancherel's Formula:}
We can relate these transforms back to our problem of computing norms of probability measures by Plancherel's formula, which states that for any $f,h: G \to \C$, we have the following expression: \\
\begin{equation} \label{plancherel}
\sum_{s\in G} f(s)h(s) = \frac{1}{|G|} \sum_{\rho \in \widehat{G}} d_\rho \text{Tr}(\widehat{f}(\rho)\widehat{h}(\rho)^*)
\end{equation}
Where $\widehat{h}(\rho)^*$ is the hermitian adjoint of $\widehat{h}(\rho)$. An accessible proof is found in Chapter 6 of Serre \cite{Serre}. We can think of the left hand side as a sort of inner product on the space of functions $G \to \C$, often denoted as $\C[G]$. In this sense, it induces a norm, which we can write as:
\[
\sum_{s\in G} f(s)^2 = \frac{1}{|G|} \sum_{\rho \in \widehat{G}} d_\rho \text{Tr}(\widehat{f}(\rho)\widehat{f}(\rho)^*)
\]
This gets us closer to translating our problem to the language of representation theory. Now, let's set $f = P - U$, where $P$ is the probability distribution in (\ref{def:prob}) and $U$ is the uniform distribution. \\
Then by (\ref{value:U}), this reduces to: \\
\[
\sum_{s\in G} [P(s) - U(s)]^2 = \frac{1}{|G|} \sum_{\rho \in \widehat{G^*}} d_\rho \text{Tr}(\widehat{P}(\rho)\widehat{P}(\rho)^*)
\]
Using the Cauchy-Schwartz inequality, we can finally write:
\begin{align*}
||P-U||^2
&= (\frac{1}{2}\sum_{s\in G} |P(s) - U(s)|)^2 \\
&= \frac{1}{4} |G|^2 (\frac{1}{|G|} \sum_{s\in G} |P(s) - U(s)|)^2 \\
&\leq \frac{1}4 |G|^2 (\sum_{s\in G} \frac{1}{|G|^2}) \sum_{s\in G} (|P(s) - U(s)|)^2 \\
&= \frac{1}{4} |G| \sum_{s\in G} (|P(s) - U(s)|)^2 \\
&= \frac{1}{4} \sum_{\rho \in \widehat{G^*}} d_\rho \text{Tr}(\widehat{P}(\rho)\widehat{P}(\rho)^*)
\end{align*}
Which is the upper bound lemma we wanted to prove.

\section{The Diaconis-Shahshahani Upper bound on convergence}
The goal of this section is to show how we use the upper bound lemma we have just shown to develop an asymptotic bound on the convergence of the shuffle in (\ref{def:prob}).
Note that from 2.2, we know that when $P$ is a class function, then $\widehat{P}(\rho) = \lambda_\rho I$ for some $\lambda_\rho \in \C$. \\
Taking traces, we get that $\lambda_\rho = \frac{1}{n} + \frac{n-1}{n} \frac{\chi_\rho(\tau)}{d_\rho}$. To follow conventions in the literature, we shall write $r(\rho):= \frac{\chi_\rho(\tau)}{d_\rho}$. This is often referred to as the \textbf{character ratio} \\
Using this fact and lemma 4, we get: \\
\begin{align}\label{mainsum}
||P^{*k}-U||^2
\leq \frac{1}{4} \sum_{\rho \in \widehat{G^*}} d_\rho \text{Tr}(\widehat{P}(\rho)\widehat{P}(\rho)^*)
= \frac{1}{4} \sum_{\rho \in \widehat{G^*}} d^2_\rho (\frac{1}{n} + \frac{n-1}{n}r(\rho))^{2k} 
\end{align}
In the last step, we used the fact that the $\chi_\rho(\tau)$ is real (matter of fact, we know that the characters of $S_n$ are all integers). Note that because $\rho(g)$ is always a unitary matrix, we have the basic inequality that $|r(\rho)| \leq 1$ and therefore, the term in the sum $(\frac{1}{n} + \frac{n-1}{n}r(\rho))^2 \leq 1$ $\forall \rho \in \widehat{G}$. This inequality is indeed tight, since upper bound is realized when $\rho$ is simply the trivial representation.\\
What's remarkable is that for most representations, $|r(\rho)|$ is actually much smaller than 1. Diaconis and Shahshahani's core idea is to partition $\widehat{G}$ into sets or "regions", in each of which we have a much better (exponential) bound on the summands in (\ref{mainsum}). To motivate their approach, we present a couple of useful facts on the structure of $\widehat{S_n}$:
\subsection{Representation theory of the symmetric group}
Since character tables are square (the number of conjugacy classes equals the number of irreducible characters), we know that $\widehat{G}$ has the same number of elements as the number of partitions of the number $n$. Turns out there's a much more direct correspondence between these two sets through the idea of \textbf{Young Tableaux}.
Chapter 7 of Diaconis' book on representation theory in probability \cite{diaconis_book} contains an introduction to the representation theory of the symmetric group. The main takeaway is that we can associate each irreducible representation of $S_n$ with a partition of the integer $n$, which we call $\lambda$. \\
To say that $\lambda$ is a \textbf{partition} of $n$, we mean that $\lambda$ is a tuple of natural numbers $(\lambda_1, \lambda_2, \ldots, \lambda_m)$ such that $\lambda_1 \geq \lambda_2 \geq \ldots \geq \lambda_m$ and that $\sum \lambda_i = n$. We write $\lambda \vdash n$ as shorthand for "$\lambda$ is a partition of the integer $n$"\\
It's often useful to visualize these representations/partitions wit the help of Young Tableaux, diagrams representing the partitions. A tableau is made up of rows of squares such that the tableau has $m$ rows, and the $i$-th row has exactly $\lambda_i$ squares. \\
The \textbf{transpose} of a partition, denoted as $\lambda'$, is a partition which corresponds to the Young Tableau obtained by "flipping" the the $\lambda$ tableau along its main diagonal. \\
For example when $n = 7$, $\lambda = (3,2,1,1)$ corresponds to the Young Tableau on the left in figure 1. Its transpose, $\lambda'$, corresponds to the one on the right: \\
\begin{figure}[h] \caption{$\lambda = (3,2,1,1)$ and its transpose}
\includegraphics[width=100mm]{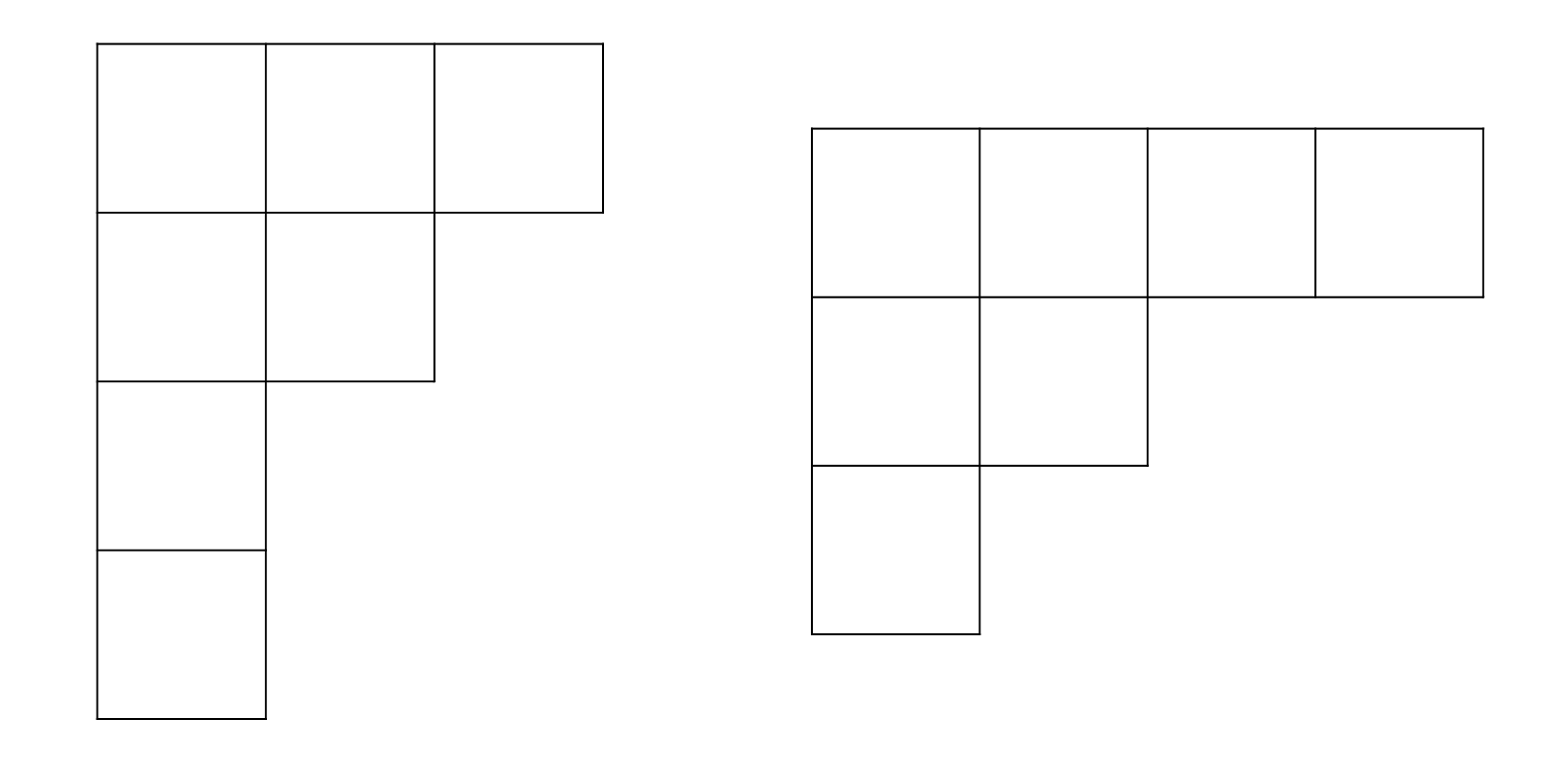}
\centering
\end{figure}
It's easy to see that when $\lambda$ has $m$ rows and its first row is of size $\lambda_1$, then $\lambda'$ has $\lambda_1$ rows and its first row is of size $m$. In our previous example, $\lambda'= (4,2,1)$
\subsection{Understanding the character ratio $r(\lambda)$}
Much of the problem of bounding the error in \ref{mainsum} comes down to estimating $r(\lambda)$ for various representations. Luckily, we actually have explicit formulas for the characters of the symmetric group. For our discussion, we're only interested in the values of the characters on the transpositions. According to \cite{James}, we have that
\begin{equation} \label{r-value}
    r(\lambda) = \frac{2}{\binom{n}{2}} \sum_{k=1}^m [\binom{\lambda_k}{2} - \binom{\lambda'_k}{2}] = \frac{1}{n(n-1)} \sum_{k=1}^m [(\lambda_j - j)^2 + (\lambda_j -j) + j(j-1)]
\end{equation}
We can induce the following partial ordering on the partitions of $n$. Let $\lambda, \mu \vdash n$. We write $\mu \trianglelefteq \lambda$ to mean that: $\lambda_1 \geq \mu_1$, $\lambda_1 + \lambda_2 \geq \mu_1 + \mu_2$ (generally, $\sum_i^k \lambda_i \geq \sum_i^k \mu_i$, $\forall 1 \leq i \leq m$). An equivalent characterization of this partial order is that $\mu \trianglelefteq \lambda$ if and only if the Young Tableau of $\lambda$ can be obtained from that of $\mu$ by iteratively picking one square from the end of a row and moving it to a row above it. \\
The following lemma will prove central to the bounds on (\ref{mainsum}):
\begin{lemma}
If $\mu \trianglelefteq \lambda$, then $r(\mu) \leq r(\lambda)$
\end{lemma}
\noindent A proof can be found in chapter 3D of Diaconis \cite{diaconis_book}, in which we induce on the operation of moving a box up in the Young Tableau to get from the Tableau of $\mu$ to that of $\lambda$ in a finite number of steps, showing that this operation increases the value of the character ratio using the formula in (\ref{r-value}). \\
The importance of this lemma comes down to the following observation. Let $\lambda^* \vdash n$ and $\lambda^* = (k,n-k)$. Then $\forall \lambda \in \widehat{G}$ such that $\lambda_1 \leq k$, we have $r(\lambda) \leq r(\lambda^*)$. Since $\lambda^*$ partitions $n$ into only 2 integers,  $r(\lambda^*)$ is quite simple to compute and acts as an upper bound on $r(\lambda)$ for a large class of elements of $\widehat{G}$. There's a balance, however, to be struck when choosing the $\lambda^*$ to work with: we want a $\lambda^*$ with a "small enough" $r(\lambda)$ as to provide a strong upper bound. On the other, we want to define a simple enough lambda with a large $\lambda_1$ as to bound a large number of summands in (\ref{mainsum}). 
\subsection{Dividing into regions}
In order to do this, Diaconis and Shahshahani resort to dividing $\widehat{G}$ into multiple "zones" in which a corresponding $\lambda^*$ generates the corresponding upper bound.\\
The partition used by Diaconis and Shahshahani is the following. Let $\lambda=(\lambda_1, \ldots, \lambda_m)$, then we can place $\lambda$ into one of the following 3 subsets of $\widehat{G}:= A_1 \cup A_2 \cup A_3$: \\
\textbf{Inner Zone $(A_1)$:} $\lambda_1 \leq \frac{n}{3}$ and $m \leq \frac{n}{3}$  \\
\textbf{Mid Zone $(A_2)$:} $\{\frac{n}{3} < \lambda_1 \leq \frac{n}{2}$ and $m \leq \frac{n}2\}$ OR $\{\frac{n}{3} < m \leq \frac{n}{2}$ and $\lambda_1 \leq \frac{n}{2}\}$ \\
\textbf{Outer Zone $(A_3)$:} $\frac{n}{2} < \lambda_1$ or $\frac{n}{2} < m$ \\
It's easy to convince oneself that each representation is an element of exactly one of these zones. The zones are illustrated in figure 2, which appears in \cite{DS}.
\begin{figure}[h] \caption{Partitioning the irreducible representations}
\includegraphics[width=62mm]{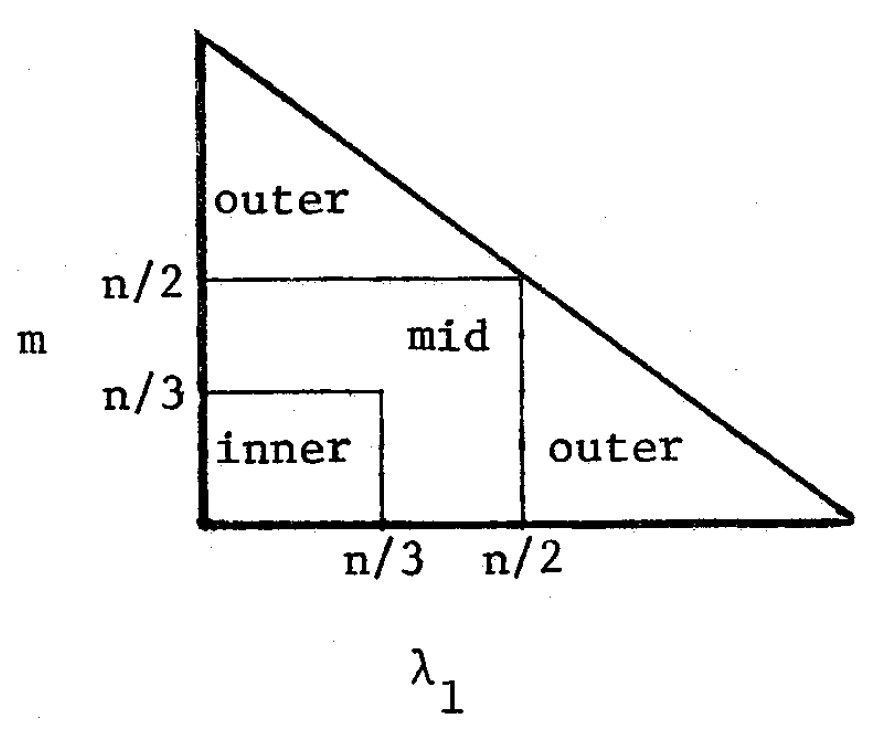}
\centering
\end{figure}
\\
\noindent The trivial representation corresponds to $\lambda = (n), m = 1$ and, while it is not part of the sum in ($\ref{mainsum}$), it appears here in the outer zone. \\
We show the following bound on the inner zone and omit the rest of the zones since our focus is on the proof techniques. The bounds of the other zones, although trickier to show, use very similar techniques.
\begin{lemma}
Let $A_1$ be defined as earlier. We have :
\begin{equation}
\sum_{\lambda \in A_1} d_\lambda^2 (\frac{1}{n} + \frac{n-1}{n} r(\lambda))^{2k} < (\frac{1}{3})^{2k} n!
\end{equation}
\end{lemma}
\noindent \textit{Proof:} \\
Fix $n$ and let $\epsilon \geq 0$ be the smallest fraction such that $b := \frac{n}{3} + \epsilon$ is an integer. Let $\lambda^*:= (b,b,n-b)$. \\
Note that $\forall \lambda \in A_1$, $\lambda \trianglelefteq \lambda^*$ and $\lambda' \trianglelefteq \lambda^*$ (since $\lambda \in A_1$ implies that $\lambda' \in A_1$). \\
Taking advantage of the symmetry of the first form of the formula for $r(\lambda)$ in (\ref{r-value}), note that we have $r(\lambda') = - r(\lambda')$ $\forall \lambda$. From this we can conclude $|r(\lambda)| \leq (\lambda^*)$ when $\lambda \in A_1$. \\
Using the second form of the formula in (\ref{r-value}), we obtain a simple closed formula for $r(\lambda^*)$:
\begin{align*}
    r(\lambda^*) &= \frac{1}{n(n-1)} [b(b-1) + (b-1)(b-2) + (n-2b-1)(n-2b-2) -2 -6] \\
    & \leq \frac{1}{(n-1)} [\frac{n}{3} - 3 + \frac{20}{3} \cdot \frac{1}{n}]
\end{align*}
Where the second line follows from expansion and algebraic multiplication, and setting $\epsilon = \frac{2}{3}$. From this we get that:
\begin{equation*}
    \frac{1}{n} + \frac{n-1}{n} r(\lambda)
    < \frac{1}{n} + \frac{n-1}{n} r(\lambda^*)
    = \frac{1}{3} - \frac{2}{n} + \frac{20}{3} \cdot \frac{1}{n^2}
    < \frac{1}{3}
\end{equation*}
Where the last inequality holds whenever $n > 3$. \\
To get a lower bound, note that since $\lambda' \in A_1$, then by the bound above:
\begin{align*}
    & r(\lambda') \leq \frac{1}{n-1}[\frac{n}{3} - 3 + \frac{20}{3} \cdot \frac{1}{n}] \\
    & r(\lambda) \geq \frac{1}{n-1}[-\frac{n}{3} + 3 - \frac{20}{3} \cdot \frac{1}{n}] \\
    & \frac{n-1}{n} r(\lambda) \geq -\frac{1}{3} + \frac{3}{n} - \frac{20}{3} \cdot \frac{1}{n^2} \\
    & \frac{1}{n} + \frac{n-1}{n} r(\lambda) \geq -\frac{1}{3} + \frac{4}{n} - \frac{20}{3} \cdot \frac{1}{n^2} \\
    & \frac{1}{n} + \frac{n-1}{n} r(\lambda) > -\frac{1}{3} \text{ for $n > 3$}
\end{align*}
From this, we can conclude that $|\frac{1}{n} + \frac{n-1}{n} r(\lambda))| < \frac{1}{3}$, which gives us the desired bound that $(\frac{1}{n} + \frac{n-1}{n} r(\lambda))^{2k} < (\frac{1}{3})^{2k}$. \\
\\
To finish the proof, we present a rudimentary bound on $d_\lambda$ that does not depend on $A_1$. We note from Corollary 5, chapter 2 in Serre \cite{Serre} that $\sum_{\rho \in \widehat{G}} d_\rho^2 = |G|$ for any group $G$. In particular, we have: $\sum_{\lambda \in A_1} d_\lambda^2 \leq \sum_{\lambda \in \widehat{S_n}} d_\lambda^2 = n!$. This gives us the desired lemma.
\subsection{Bounding the error terms}
\noindent This technique can be applied to the rest of the regions. We face the need to further subdivide $A_3$ into 3 distinct regions, which we call here $B_1, B_2,$ and $B_3$. We end up getting the following bounds for summand corresponding to $\lambda$:
\begin{align*}
    A_1: & (\frac{1}{3})^{2k}n! \\
    A_2: & \exp{(\pi (\sqrt{\frac{2}{3}n)})}4^n(\frac{1}{2})^{2k}n^\frac{2n}{3} \\
    B_1: & \exp{(\pi (\sqrt{\frac{2}{3}n)})}4^n(\frac{n}{2})!(\frac{29}{50} + \frac{1}{5n} + \frac{4}{n^2})^{2k} \\
    B_2: & \exp{(-\frac{4k}{n})} \sum_{j=0}^{\frac{3}{10}n} \frac{p(j)}{j!}\exp{2j^2\frac{log(n)}{n}} \\
    B_3: & e^{-4c} \sum_{j=1}^{\frac{3}{10}n} \frac{p(j)}{j!}\exp{(2j(j-1)\frac{log(n)}{n})} 
\end{align*}
Here, $c := \frac{k}{n} - \frac{1}{2}\log(n)$ and $p(n)$ is the number of partitions of $n$, which is precisely $|\widehat{G}|$ in our case. \\
In order to finish the proof for theorem \ref{thm:DS}, we would need to simply show that the sum of these four terms can be bounded from above by $ae^{-4c}$ for some constant $a$ which does not depend on $c$ or $n$. \\
Just as a reminder, once we accomplish this, then 
Once again, we show how this is done for $A_1$ for illustration's sake. The key is to use Stirling's approximation to bound factorials by exponentials. Stirling's approximation tells us that
\[
    \log(n!) = n \log(n) - n + O(\log(n))
\]
From this, we get the following "coarser" inequality, which suffices, which is that:
\[
    n! \leq e^{n \log n}
\]
Using this, we can bound the $A_1$ term:
\begin{align*}
    (\frac{1}{3})^{2k}n! 
    = (\frac{1}{3})^{n \log n + 2cn}n!
    \leq (\frac{1}{3})^{n \log n + 2cn} e^{n \log n }
    = (\frac{e}{3})^{n \log n} (\frac{1}{3})^{2cn}
    < (\frac{1}{3})^{2cn}
    < e^{-2cn}
    < e^{-4c}
\end{align*}
When $n > 3$.
The bounds on the rest of the region, albeit a bit more technical, are in the same spirit of this. One additonal tool which comes up is the classic asymptotic expression for $p(n)$, obtained by Ramanujan and Hardy, which is that:
\[
p(n) \sim \frac{1}{4n\sqrt{3}} \exp{(\pi \sqrt{\frac{2n}{3}})}
\text{ as} \hspace{10pt} n \to \infty
\]

\section{Determining the limiting profile}
In \cite{Teyssier} Teyssier determines the limiting profile of the total variation distance in (\ref{eq:2}). We focus here on a couple of notable techniques in which this is done.\\
As a reminder, the statement we aim to show is the following, using the same definition for $P$ as in (\ref{def:prob}), fixing $c$ and $k(n,c)$ as previously, we have:
\[
\underset{n \to \infty}{\lim} ||P^{*k} - U||_{TV} = ||\text{Poiss}(1+e^{-2c}) - \text{Poiss}(1)||_{TV}
\]
While we will not present the proof in its entirety, we will focus on two key ideas from the proof which might be of relevance to showing limiting profiles for other problems which center around the convergence of a random walk to its stationary distribution. \\
Diaconis aand Shahshahani note in \cite{DS} that the representations with the largest contributions to the sum in (\ref{mainsum}) are the ones with relatively large $\lambda_1$. One of Teyssier's contributions is making this notion more precise by proving a bound - that is uniform in $n$ - on the contribution of representations which have a bounded $\lambda_1$. This will be detailed in Lemma 9. \\
The second notable contribution he makes is a modification of the Diaconis-Shahshahani upper bound lemma proven in section 2 of this paper, which is the following statement:
\begin{lemma} \label{modifiedULM}
Let $G$ be any group, $P$ be a class function, $s_\rho := \frac{\text{Tr}(\widehat{P}(\rho))}{d_\rho}$ $\forall \rho \in \widehat{G}$, and $S \subseteq \widehat{G}^*$then:
\begin{equation}
    \left| ||P^{*k} - U|| - \frac{1}{2|G|}\sum_{g \in G} \left |\sum_{\rho \in S} d_\rho s_\rho^t \overline{\chi_\rho(g)} \right |   \right|
    \leq \frac{1}{2} \sum_{\rho \in \widehat{G}^* \slash S} d_\rho|s_\rho^t|
\end{equation}
\end{lemma}

The proof for this is rather straightforward and relies on the Fourier inversion formula, which can be seen as a special case of Plancherel's Formula:
\begin{lemma} (Fourier inversion formula)
Let $f \in \C[G]$, then $\forall g \in G$
\begin{equation}
    f(g) = \sum_{\rho \in \widehat{G}} \frac{d_\rho}{|G|} \text{Tr}(\rho(g)^*\widehat{f}(\rho))
\end{equation}
\end{lemma}
\noindent \textit{Proof (Lemma 8):} Fix $g \in G$ and let $\delta_g \in \C[G]$ be the kronecker delta function, then for $\rho \in \widehat{G}$:
\[
\widehat{\delta_g}(\rho) = \sum_{s \in G} \delta_g(s) \rho(s) = \rho(g)
\]
So using (\ref{plancherel}), we get:
\[
f(g) = \sum_{s\in G} f(s)\delta_g(s)
= \sum_{\rho \in \widehat{G}} \frac{d_\rho}{|G|} \text{Tr}(\widehat{\delta_g}(\rho)^*\widehat{f}(\rho))
= \sum_{\rho \in \widehat{G}} \frac{d_\rho}{|G|} \text{Tr}(\rho(g)^*\widehat{f}(\rho))
\]
We can now go back to proving lemma 7, which is rather straightforward albeit a bit notation heavy. We begin by expressing the first term in the difference as a sum over $\widehat{G}^*$ using the Fourier inversion formula:
\begin{align*}
||P^{*k} - U||
& = \frac{1}{2} \sum_{g \in G} |P^{*k}(s) - U(s)|
= \frac{1}{2} \sum_{g \in G} |P^{*k}(s) - U(s)|
= \frac{1}{2} \sum_{g \in G} |\sum_{\rho \in \widehat{G}} \frac{d_\rho}{|G|} \text{Tr}(\rho(g)^*(\widehat{P^{*k} - U})(\rho))| \\
& = \frac{1}{2} \sum_{g \in G} |\sum_{\rho \in \widehat{G}^*} \frac{d_\rho}{|G|} \text{Tr}(\rho(g)^*(\widehat{P^{*k}}(\rho))|
= \frac{1}{2} \sum_{g \in G} |\sum_{\rho \in \widehat{G}^*} \frac{d_\rho}{|G|} \text{Tr}(\rho(g)^*\widehat{P^{*k}}(\rho)| \\
& = \frac{1}{2|G|} \sum_{g \in G} |\sum_{\rho \in \widehat{G}^*} d_\rho s_\rho^t \overline{\chi_\rho(g)}| 
\end{align*}
Using this, we can see that:
\begin{align*}
    & \left| ||P^{*k} - U|| - \frac{1}{|G|}\sum_{g \in G} \left |\sum_{\rho \in S} d_\rho s_\rho^t \overline{\chi_\rho(g)} \right |   \right|
    = \frac{1}{2|G|} \left| \sum_{g \in G} |\sum_{\rho \in \widehat{G}^*} d_\rho s_\rho^t \overline{\chi_\rho(g)}| - |\sum_{\rho \in S} d_\rho s_\rho^t \overline{\chi_\rho(g)}|  \right| \\
    & \leq \frac{1}{2|G|} \sum_{g \in G}  \left| |\sum_{\rho \in \widehat{G}^*} d_\rho s_\rho^t \overline{\chi_\rho(g)}| - |\sum_{\rho \in S} d_\rho s_\rho^t \overline{\chi_\rho(g)}|  \right|
    \leq \frac{1}{2|G|} \sum_{g \in G} |\sum_{\rho \in \widehat{G}^* \slash S} d_\rho s_\rho^t \overline{\chi_\rho(g)}| \\
    & \leq \frac{1}{2|G|} \sum_{g \in G} \sum_{\rho \in \widehat{G}^* \slash S} |d_\rho s_\rho^t \overline{\chi_\rho(g)}|
    = \frac{1}{2} \sum_{\rho \in \widehat{G}^* \slash S} d_\rho |s_\rho^t| \sum_{g \in G} \frac{1}{|G|} |\overline{\chi_\rho(g)}| \\
    & \leq \frac{1}{2} \sum_{\rho \in \widehat{G}^* \slash S} d_\rho |s_\rho^t| \sqrt{\sum_{g \in G} \frac{1}{|G|^2}} \cdot \sqrt{\sum_{g \in G} |\overline{\chi_\rho(g)}|^2}
    = \frac{1}{2} \sum_{\rho \in \widehat{G}^* \slash S} d_\rho |s_\rho^t|
\end{align*}
Where we repeatedly used the triangle inequality. The last inequality is indeed where we use Cauchy-Schwartz. We can now use this lemma to show that a large subset of $\widehat{G}$ contributes a negligible mass, in the followng sense:
\begin{lemma}
Let $P$ be defined as in (\ref{def:prob}), then $\forall
\epsilon > 0$, $c \in R$, $\exists M \geq 1$, $n_0$, s.t $n \geq n_0$ implies that:
\[
\sum_{\lambda \in \widehat{G}^* \slash S_M} d_\rho |s_\rho^t| \leq \epsilon
\]
Where $S_M:= \{\lambda \in \widehat{G} : \lambda_1 \geq n-M\}$
\end{lemma}
\noindent The proof of this lemma proceeds similarly to that of \ref{thm:DS} so we omit it here. In essence, the set $\widehat{G}^* \slash S_M$ is divided into zones, each of which is bounded above by some function of $\epsilon$ and $n$, which is ultimately $o(1)$ in $n$. \\
To see why this result is significant for our problem, note that $\widehat{G}^* \slash S_M$ contains the vast majority of the irreducible representation. More concretely, we have the following lemma, which we shall now prove:
\begin{lemma}
Let $M$ and $S_M$ be defined as previously, then $\underset{n \to \infty}{\lim} \frac{|S_M|}{|\widehat{G}|} = 0$
\end{lemma}
\noindent \textit{Proof:}
Let $\lambda \in S_M$. By the 1-to-1 correspondence of the irreducible representations of $S_n$ and its conjugacy classes, we can identify $\lambda$ with a certain cycle shape in $S_n$. Let $\lambda_1 = k$. For sufficiently large $n$ ($n > 2M$), $\lambda$ must have exactly 1 cycle of length at least $n-M$. The number of unique conjugacy classes in $S_n$ with exactly 1 $k$-cycle is simply $p(n-k)$, where $p$ is the partition function. From this we get that:
\[
|S_M| = \sum_{k=n-M}^n p(n-k) = \sum_{k = 0}^M p(k) \leq Mp(M)
\]
Where the last inequality follows from the fact that $p$ is strictly increasing. Therefore, we get that:
\[
\frac{|S_M|}{|\widehat{G}|} \leq \frac{Mp(M)}{p(n)}
\]
which clearly decays to 0 exponentially quickly. To be a bit precise, Hardy and Ramanujan present simple lower and upper bounds for the partition function in section 2 of \cite{HR}: For universal constants $A,B \in \R$, we have that $e^{A\sqrt{n}} \leq p(n) \leq e^{B\sqrt{n}}$. Applied to our problem, this mean:
\[
\log \left( \frac{|S_M|}{|\widehat{G}|} \right) < B\sqrt{m} -  A\sqrt{n}
\]
\section*{Concluding Remarks}
The techniques used in \cite{Teyssier} could be relevant to determining the limiting profile of other random walks. In particular, the modified Diaconis-Shahshahani upper bound in (\ref{modifiedULM}) has a fairly wide scope, since it makes no stipulations about the nature of the finite group $G$ or the function $P \in \C[G]$ (beyond the fact that it is a class function). \\
In a recent paper \cite{Nestoridi}, Nestoridi and Olesker-Taylor use this upper bound to derive the a couple of limiting profiles, including that of the k-cycle random walk on $S_n$ and, in a sense, generalizing Teyssier's work (the walk considered in \cite{Nestoridi} is albeit different in one central way, which is that no mass is placed on the identity element of $S_n$ in defining the distribution, which leads to the need for making certain considerations around parity, which we didn't have to make when studying the distribution in (\ref{def:prob}). \\
We hope that these techniques can be applied to other problems. As of the time of writing, the limiting profile of the random walk on the special linear group $SL_n(F_q)$ (here, $q$ is prime) has yet to be determined. The cutoff phenomenon, however, has been studied by Hildebrand in \cite{hildebrand_1992}. In his proof for the cutoff phenomenon, he employs a zone division argument similar to the one used in Section 3.4. One way the cutoff profile problem could be attacked is by carefully examining the bounds and determining a suitable partition into a set corresponding to an "error term" (the set of representations whose corresponding summands have contribution the TV distance that is uniformly bounded by $\epsilon$ in the limit as $n$ tends to infinity), and a "limiting term" (the irreducible which contribute the most, whose corresponding summands can be manipulated and shown to converge to the limiting profile).

\newpage
\section*{Acknowledgements}
I would like to thank my Advisor, Prof. Evita Nestoridi, for introducing me to card shuffles and the work of Prof. Persi Diaconis on the representation theory of the symmetric group and for her guidance throughout this project.

\newpage
\printbibliography 
\end{document}